\documentclass[12pt]{amsart}
\usepackage{amsfonts, amsmath, amsthm}
\usepackage{psfig}

\newtheorem{pro}{Proposition}[section]
\newtheorem{thm}[pro]{Theorem}
\newtheorem{lem}[pro]{Lemma}

\newtheorem{cor}[pro]{Corollary}

\theoremstyle{definition}

\theoremstyle{remark}

\font\cal=cmsy10
\begin{document}

\title{Length multiplicities of hyperbolic 3-manifolds}
\author{Joseph D. Masters}

\begin{abstract}
Let $M = \mathbb{H}^3/\Gamma$ be a hyperbolic 3-manifold,
 where $\Gamma$ is a non-elementary
 Kleinian group.  It is shown that the length spectrum of $M$
 is of unbounded multiplicity.
\end{abstract}

\maketitle

\section{Introduction}
 This paper is concerned with length multiplicities in hyperbolic
 3-manifolds, or more generally, in hyperbolic 3-orbifolds.
 Let $M = \mathbb{H}^3/\Gamma$ 
  be a hyperbolic 3-orbifold, where $\Gamma$ is a non-elementary
 Kleinian group. We say that $\gamma \in \Gamma$ is
 \textit{loxodromic} if $tr^2 \gamma \not\in [0,4]$ (note that this
 includes ``hyperbolic'' elements).
 Every loxodromic element $\gamma \in \Gamma$ has an associated 
 \textit{complex length}, denoted $\ell_0(\gamma) = \ell + i\theta$,
 which describes the action of $\gamma$ on $\mathbb{H}^3$:
 along its invariant axis $\gamma$ translates a distance $\ell$
 and rotates an angle $\theta$.  We say that a complex length
 has \textit{multiplicity n} if it is shared by exactly $n$
 conjugacy classes in $\Gamma$.
  We define the \textit{complex length spectrum}, $\cal{L}(M)$, to be the
 set of complex lengths of loxodromic elements of $\Gamma$,
 counted with multiplicity.

  The \textit{real length spectrum} of $M$ is defined
 to be the set of lengths of closed geodesics of $M$.
 The real length spectrum of $M$ is essentially the
 real part of $\cal{L}(M)$-- the only difference being that the classes
 of $\gamma$ and $\gamma^{-1}$ are now equivalent.

  Since $\Gamma \subset PSL_2(\mathbb{C})$,  $tr(\gamma)$ is
 well-defined, up to sign, for any $\gamma \in \Gamma$.
  The connection between
 traces and lengths is given by the following formula:

\begin{equation}
\ell_0(\gamma) = 2cosh^{-1}(\frac{tr \,\, \gamma}{2}).
\end{equation}

 Following [GR], we define, for any group $\Gamma$, 
 the \textit{trace class} of an element
 $\gamma \in \Gamma$ to be the set of elements $\gamma^{\prime} \in \Gamma$
 for which $tr \rho(\gamma^{\prime}) = tr \rho(\gamma)$ for \textit{all}
 representations $\rho : \Gamma \rightarrow SL_2(\mathbb{C})$.
  We define the \textit{stable multiplicity} of $\gamma$
 to be the number of conjugacy classes in the trace class of $\gamma$. 
 Recall that a representation into $SL_2(\mathbb{C})$ is
 called \textit{irreducible} if its image fixes no 1-dimensional
 subspaces of $\mathbb{C}^2$.

 We shall prove:

\begin{thm} \label{trace}
Let $\Gamma$ be a
 finitely generated group which admits an infinite
 irreducible representation into $SL_2(\mathbb{C})$.  Then 
 $\Gamma$ has trace classes of unbounded stable multiplicity.
\end{thm}

In Section 2 we prove that Theorem \ref{trace} has the following consequence:

\begin{thm} \label{intro}
Let $M = \mathbb{H}^3/\Gamma$ be a hyperbolic 3-orbifold, where
 $\Gamma$ is a finitely generated, non-elementary Kleinian group.
Then $\cal{L}(M)$ is of unbounded multiplicity.
\end{thm}

\begin{cor} If $M$ is a finite-volume, complete hyperbolic \newline
 3-manifold, then $\cal{L}(M)$ is of unbounded multiplicity.
\end{cor}

\begin{cor} If $M$ is a finite-volume, complete hyperbolic\\
 3-manifold, then the real length spectrum of $M$ is of unbounded multiplicity.
\end{cor}

 The analogous statement for hyperbolic surfaces was proved by
 Randol (see [R]). We have reproduced the short proof in Section 3.

 It is well-known that the length spectrum of an
 arithmetic hyperbolic 3-manifold has unbounded multiplicity.
  In fact, if $M$ is arithmetic, then the mean multiplicity,
 $n(\ell_0)$, of a length grows exponentially with $\ell_0$
 (see [M]).
 
  Recent interest in length multiplicities
 of hyperbolic 3-manifolds has been sparked by connections
 with chaotic quantum systems.  See [Sar] for more information.

I would like to thank Professors Fernando Rodriguez Villegas,
 John Tate and Cameron Gordon for useful conversations, Professor
 Peter Sarnak for valuable correspondence,
 the referee for helpful comments
 and Professor
 Alan Reid for his help and encouragement.

\section{Theorem \ref{trace} implies Theorem \ref{intro}}

In this section, we prove:\\
\\
\textit{Claim:} Theorem \ref{trace} implies Theorem \ref{intro}\\
\\
\textit{Proof of claim:}  By Selberg's Lemma, $\Gamma$ has a torsion-free
 subgroup $\Gamma^{\prime}$ of finite index $n$, say.
 Any set of more than $n$ elements which are pairwise non-conjugate in
 $\Gamma^{\prime}$
 must contain at least two elements which are 
 non-conjugate in $\Gamma$.  Also, any representation of $\Gamma$
 restricts to a representation of $\Gamma^{\prime}$.  Therefore the 
 stable multiplicities in $\Gamma^{\prime}$
 will increase by at most a factor of $n$.  So it
 is enough to consider the case where $\Gamma$ is torsion-free.

Let $\pi :SL_2( \mathbb{C}) \rightarrow PSL_2(\mathbb{C})$
 be the natural projection. 
 By [T], a non-elementary, torsion-free
 Kleinian group $\Gamma$ admits a faithful
 representation $\rho :\Gamma \rightarrow \pi^{-1}(\Gamma)$
 such that $\pi \rho = id$.  Note that
 $\rho$ preserves trace (up to sign), and that
 $\rho$ is irreducible, since $\Gamma$ is non-elementary. Then, in the
 case that $\Gamma$ contains no parabolic elements,
 Theorem \ref{intro} is now an easy consequence of Theorem \ref{trace}
 and Equation 1.

 In general, however, we must make sure that $\Gamma$ has trace classes
 with an unbounded
 number of \textit{loxodromic} conjugacy classes.
 The claim will be proved once
 we show that any trace class in $\Gamma$ can
 contain only a bounded number of conjugacy classes of parabolic elements.
 This is done in the following lemma:

\begin{lem} \label{Klein}
Let $\Gamma$ be a finitely generated Kleinian group. 
 Then there is
 an integer $N > 0$ such that no trace class of $\Gamma$ contains
 more than $N$ conjugacy classes of parabolic elements.
 Moreover, if $\Gamma$ is geometrically finite, then a trace
 class of $\Gamma$ can contain at most two conjugacy class of
 parabolic elements.
\end{lem}

The proof of this lemma will require some definitions.

Let $\Gamma$ be a finitely generated group.  The \textit{space of
 characters},
 $V(\Gamma)$, is the set of all characters of representations of
 $\Gamma$ into $SL_2(\mathbb{C})$.
  By [CS], $V(\Gamma)$ has the structure of an affine algebraic set
 defined over $\mathbb{Q}$.
 The character of a representation $\rho$ is denoted $\chi_{\rho}$.
  
 A \textit{quasi-conformal deformation of} $\Gamma$ is a
 representation of $\Gamma$ into $PSL_2(\mathbb{C})$
 which is induced by a quasiconformal homeomorphism of the
 Riemann sphere $\hat{\mathbb{C}}$.  We say that a represntation
 of $\Gamma$ into $SL_2(\mathbb{C})$ is quasiconformal if it
 is the lift of a quasiconformal deformation of $\Gamma$ into 
 $PSL_2(\mathbb{C})$.

\begin{proof} (of Lemma \ref{Klein})
 Using Selberg's Lemma as above, we may assume $\Gamma$ is
 torsion-free, so the identity representation lifts to a representation
 $\rho_0 :\Gamma \rightarrow SL_2(\mathbb{C})$.

  By the compact core
 theorem ([Sc]), $\Gamma$ can contain only finitely many conjugacy
 classes of maximal parabolic subgroups.
 Let\\
 $\alpha_1, \beta_1, ..., \alpha_n, \beta_n$
 generate the conjugacy classes of rank 2 maximal parabolic subgroups
 and $\gamma_1, ..., \gamma_m$ generate the conjugacy classes
 of rank 1 maximal parabolic subgroups.

 First, suppose $\Gamma$ is finite covolume,
 so there are no $\gamma_i$'s.  We shall handle this case
 with Thurston's hyperbolic Dehn surgery theory.
 
 Let $V_0(\Gamma)$ denote the irreducible component of $V(\Gamma)$
 containing $\chi_{\rho_0}$.
 Consider the map
 $\tau: V_0(\Gamma) \rightarrow \mathbb{C}^n$ defined by
 $$\tau( \chi_{\rho})
 =  (\chi_{\rho}(\alpha_1), \chi_{\rho}(\alpha_2),...,\chi_{\rho}(\alpha_n))
 =  (tr \rho(\alpha_1), tr \rho(\alpha_2), ... , tr \rho(\alpha_n)).$$
 By Chapter 5 of [T], the image of $\tau$ covers an open neighborhood
 $U$ of $(2, 2, ..., 2)$.  
 So given two parabolic elements on distinct cusps,
 we can make one loxodromic while the other remains parabolic. 
 In particular, for any $i \neq j$, and any integers
 $m_1,n_1,m_2,n_2$ (not all 0) there is a representation $\rho$
 for which
 $tr \rho(\alpha_i^{m_1} \beta_i^{n_1})
  \neq tr \rho(\alpha_j^{m_2}\beta_j^{n_2})$,
 so $\alpha_i^{m_1}\beta_i^{n_1}$
 and $\alpha_j^{m_2}\beta_j^{n_2}$ are not in the same trace class.
 
 Now suppose we are given two distinct, non-trivial parabolic
 elements $\alpha_i^{m_1}\beta_i^{n_1}$
 and $\alpha_i^{m_2}\beta_i^{n_2}$ on the same cusp. 
 Suppose also that the elements are not in the same
 cyclic subgroup.
 Then Thurston's hyperbolic Dehn surgery theory
 shows that there is a representation taking
 one of the elements to the identity and the other to a loxodromic
 element, so they are in distinct trace classes.
 If the two elements are in the same cyclic subgroup,
 then by mapping them to loxodromics we see that they are
 in distinct trace classes, provided they are not inverses
 of each other.

It follows that
 we can find characters of representations in $\tau^{-1}(U)$ which differ
 on any distinct pair of parabolic elements in $\Gamma$ which
 are not inverses of each other,
 concluding the case where $\Gamma$ has finite covolume.

 If $\Gamma$ is geometrically finite, then by [Br]
 there is a quasi-conformal deformation $\Gamma^{\prime}$ of $\Gamma$
 and a finite covolume Kleinian group $\Gamma^*$
 for which $\Gamma^{\prime} \subset \Gamma^*$.
 Since a quasi-conformal deformation takes parabolics to parabolics,
 the proof of this case now follows from the proof
 of the finite covolume case.

 In general, if $\Gamma$ is any finitely generated Kleinian group,
 then $\Gamma$ has a faithful discrete representation $\rho$
 for which $\rho(\Gamma)$ is geometrically finite,
 $\rho(\alpha_i)$ and $\rho(\beta_i)$ are parabolic for all $i$,
 and $\rho(\gamma_i)$ is loxodromic for all $i$;
 this is Thm 2.3 of [A] and follows from the Scott core theorem
 and the Thurston uniformization theorem.
 By the geometrically finite case, any two elements of the form
 $\alpha_i^r \beta_j^s$ in the same trace class of $\Gamma$
 must be inverses of each other.  So a trace class
 of $\Gamma$ can contain at most two elements of the form
 $\alpha_i^r\beta_j^s$.
  And since $\rho(\gamma_i)$ is loxodromic,
 $tr \rho(\gamma_i^r) \neq tr \rho(\gamma_i^s)$ if $|r| \neq |s|$,
 so $\gamma_i^r$ and $\gamma_i^s$ are not in the same trace class in $\Gamma$,
 and a trace class in $\Gamma$ can contain at most $2m$
 elements of the form $\gamma_i^r$ (as $i$ goes from 1 to $m$).
  Therefore a trace class in $\Gamma$ can contain at most
 $2(m+1)$ conjugacy classes of parabolic elements, and the lemma is proved.
\end{proof}

\section{Background and the idea of the proof}
 
 Suppose $\rho$ is an irreducible representation of $\Gamma$ into
 $SL_2(\mathbb{C})$.
 We first must find elements in $\rho(\Gamma)$ with the same trace.
 This can be done as in [H], with the aid of simple trace identities.
 For example, it is proved in [H] that
 $tr(a^2bab^{-1}) = tr(ba^2b^{-1}a)$ for any
 elements $a$ and $b$ in $SL_2(K)$, where $K$ is any field.
 Then $a^2bab^{-1}$ and $ba^2b^{-1}a$ are in the same trace class
 in $\Gamma$. In Section 4, we will use these identities
 to construct sequences of words in $\Gamma$, 
 all in the same trace class.

 The problem, then, is to show that
 these words are not conjugate in $\Gamma$.
 This is done by finding a homomorphism of $\Gamma$
 onto a finite group $G$, and showing that the images of the words
 are not conjugate in $G$.

 This technique is nicely illustrated in the 2-dimensional case.
 The proof we give of the following theorem is a slight modification
 of the one which appears in [R].

\begin{thm} \label{2d} \textnormal{(Randol)}
Let $M = \mathbb{H}^2/\Gamma$ be a finite-volume hyperbolic surface.
 Then $\Gamma$ contains trace classes of unbounded stable multiplicity,
 and $\cal{L}(M)$ contains lengths of unbounded multiplicity.
\end{thm}

\begin{proof}

First, let us assume $M$ is compact. 

 $\Gamma \subset PSL_2(\mathbb{R})$ is a Fuchsian group, which can be
 embedded into $SL_2(\mathbb{R})$ so that traces are preserved up to sign.
 Then by Equation 1,
  it is enough to prove that $\Gamma \subset SL_2(\mathbb{R})$
 has trace classes of unbounded stable multiplicity.

 $\Gamma$ has the standard presentation:

\begin{equation*}
\Gamma = <a_1,b_1,a_2,b_2,...,a_g,b_g \, | \, (a_1b_1a_1^{-1}b_1^{-1})...(a_gb_ga_g^{-1}b_g^{-1}) = 1>.
\end{equation*}

Note that there is a natural surjection
\begin{equation*}
\phi:\Gamma \rightarrow <a_1> * <a_2> *... *<a_g>.
\end{equation*}

In particular, $a_1$ and $a_2$ generate a free group $F$.

It follows from [H] that if $x$ and $y$ freely generate a free subgroup $F$
 of $\Gamma$, then for any N there are words
 $w_1, ..., w_N$ in $x$ and $y$
 such that\\
\\
1. $w_i(x,y)$ and $w_j(x,y)$ are in the same trace class
   $ \, \forall \, i,j \leq N$\\
2. $w_i(x,y)$ is not conjugate in $F$ 
   to $w_j(x,y) \, \forall \, i \neq j$.\\

Consider the words $w_1(a_1, a_2), ..., w_N(a_1,a_2)$.
By 1,  we have\\
 $w_i(a_1, a_2)$ and $w_j(a_1,a_2)$ are in
 the same trace class $\forall i,j$.
By 2, we have that any distinct pair $w_i(a_1,a_2)$,
$w_j(a_1,a_2)$ are not conjugate in $F$, hence their images are not conjugate
 in $<a_1> * ... *<a_g>$, hence $w_i(a_1,a_2)$ and $w_j(a_1,a_2)$
 are not conjugate in $\Gamma$.
 This proves the theorem in the compact case.

 If $M$ is non-compact, then $\Gamma$ is free.  Then the same proof
 works to show that $\Gamma$ contains trace classes of unbounded stable
 multiplicity.  The only complication is in passing to the statement
 about $\cal{L}(M)$; for the elements of the trace
 class may be parabolic.  
 However, $\Gamma$ admits faithful
 representations into $SL_2(\mathbb{C})$ for which every element
 becomes loxodromic, and therefore, as in the proof of Lemma \ref{Klein},
 we see that a trace class can contain at most $2n$ conjugacy classes
 of parabolic elements, where $n$ is the number of cusps,
 and the result follows.

\end{proof}

The proof in three dimensions is more complicated, as
 hyperbolic 3-manifold groups do not generally surject onto non-abelian
 free groups.  For example, if the manifold has zero first
 Betti number, such as a non-zero surgery on the figure-eight knot, then its
 fundamental group cannot surject onto any free group.

However, hyperbolic 3-manifold groups do surject onto groups
 of the form $PSL_2(\mathbb{F}_{p^i})$,
 where $\mathbb{F}_{p^i}$ denotes the finite
 field of order $p^i$.  In Sections 5 and 6,
 we shall review the construction of these homomorphisms.

So the idea is to use [H] to construct
 a sequence of words $w_i(a,b)$ in the free group $F$ on $a$ and $b$
 which are not conjugate in $F$ but which are in the same trace class of $F$.
 We map these words into $\Gamma$; their images will be in
 the same trace class of $\Gamma$.  Then we map the words from $\Gamma$ into
 a group of the form
 $PSL_2(\mathbb{F}_p)$, and hope that
 these images will be non-conjugate in $PSL_2(\mathbb{F}_p)$,
 so the words will be non-conjugate in $\Gamma$.
  However, by [H], the traces of the images are
 equal (up to sign) in $PSL_2(\mathbb{F}_p)$,
 and it is nearly true that two elements of $PSL_2(\mathbb{F}_p)$
 are conjugate if and only if they have
 the same trace.  Therefore
 care is needed in the choice of the words and the primes.

\section{Trace identities}

All of the trace identities which we shall use are ultimately
 based on the following lemma.

\begin{lem}\textnormal{(Horowitz)} \label{Hor}
Suppose $K$ is a field, $a,b \in SL_2(K)$ with $tr \, a = tr\, b$,
  and $W(x,y)$ is a word in $x$ and $y$.  Then\\
 $tr(W(a,b))=tr(W(b,a))$.
\end{lem}

\begin{proof}
It is proved in [H] (see also [CS]) that there is a
 3-variable polynomial $P$ over $K$ such that
 $tr(W(a,b)) = P(tr \, a, tr \, b, tr \, ab)$, for any elements
 $a,b \in SL_2(K)$.
 Then $tr(W(b,a)) = P(tr \, b, tr \, a, tr \, ba)$. 
 We have assumed that $tr \, a= tr \, b$,
 and for any matrices in $SL_2(K)$, $tr \, ab = tr\, ba$,
 so the lemma follows.
\end{proof}

In Section 3, we gave the example $tr(a^2bab^{-1}) = tr(ba^2b^{-1}a)$
 for any elements $a,b \in SL_2(K)$.  This follows
 by setting $W(x,y) = x^2y$ and noting that
 $tr(W(a,bab^{-1}))= tr(W(bab^{-1},a))$
 by Lemma \ref{Hor}.

We shall now construct the required elements of the same trace
 in $\rho(\Gamma)$ which we need to prove Theorem \ref{trace}.
  We remark that one can construct much simpler sequences of
 elements of equal trace; however, in Section 6 we shall require
 the words to be of this special form in order to prove they are
 non-conjugate.\\
\\
  We now recursively define words
 $w_{n,i}$, for $i \leq n+1$. 
 In what follows, we routinely supress the dependence
 of these words on $a$,$b$, $p_i$, $q_i$ and $k_i$;
 we will explain how to choose them later.\\
\\
Let
$W_n(x,y) = (x^{p_n-1+q_n}y^{-q_n})^{k_n}x(x^{p_n-1+q_n}y^{-q_n})x^{-1}$\\
\\
$\bar{W}_n(x,y)= x(x^{p_n-1+q_n}y^{-q_n})^{k_n}x^{-1}(x^{p_n-1+q_n}y^{-q_n})$\\
\\
\\
$w_{1,1} = W_1(a,b)$\\
\indent \indent $ =  (a^{p_1-1+q_1}b^{-q_1})^{k_1}a(a^{p_1-1+q_1}b^{-q_1})a^{-1}$\\
\\
$w_{1,2} = \bar{W}_1(a,b)$.\\
\indent \indent $= a(a^{p_1-1+q_1}b^{-q_1})^{k_1}a^{-1}(a^{p_1-1+q_1}b^{-q_1})$\\
\\
\\
$w_{2,1} = W_2(w_{1,1}, w_{1,2})$\\
\indent \indent $=  (w_{1,1}^{p_2 -1+ q_2}w_{1,2}^{-q_2})^{k_2}$
          $w_{1,1}(w_{1,1}^{p_2-1+q_2}w_{1,2}^{-q_2})w_{1,1}^{-1}$\\
\\
$w_{2,2} = \bar{W}_2(w_{1,1}, w_{1,2})$\\
\indent \indent $= w_{1,1}(w_{1,1}^{p_2 -1+ q_2}w_{1,2}^{-q_2})^{k_2}w_{1,1}^{-1}$
      $(w_{1,1}^{p_2-1+q_2}w_{1,2}^{-q_2})$\\
\\
$w_{2,3} = W_2(w_{1,2},w_{1,1})$\\
\indent \indent $= (w_{1,2}^{p_2-1+q_2}w_{1,1}^{-q_2})^{k_2}$
          $w_{1,2}(w_{1,2}^{p_2-1+q_2}w_{1,1}^{-q_2})w_{1,2}^{-1}$\\        
\\
 Assume that $w_{n-1,i}$ has been defined for $i \leq n$,
 and that $w_{n-1,1}$ and $w_{n-1,2}$ are both words in
 $w_{n+1-i,1}$ and $w_{n+1-i,2}$ for each $i$ with $3 \leq i \leq n$
 (note that this property is vacuous for $n = 2$, the base case
 of the recursion). 
\\
\\
Define
$w_{n,1} = W_n(w_{n-1,1}, w_{n-1,2})$\\
\\
$w_{n,2} = \bar{W}_n(w_{n-1,1},w_{n-1,2})$\\
\\
\\
We claim that $w_{n,1}$ and $w_{n,2}$
 are both words in $w_{n+2-i,1}$ and $w_{n+2-i,2}$
 for $3 \leq i \leq n+1$.

 Indeed, since $w_{n-1,1}$ and $w_{n-1,2}$ are words in $w_{n+1-i,1}$
 and $w_{n+1-i,2}$ for $3 \leq i \leq n$, then
 $w_{n,1}$ and $w_{n,2}$ are both words in $w_{n+2-i,1}$
 and $w_{n+2-i,2}$ for $4 \leq i \leq n+1$. 
 And for $i = 3$, it is obvious that
 $w_{n,1}$ and $w_{n,2}$ are words in $w_{n-1,1}$ and $w_{n-1,2}$.\\
\\
Then we define:\\
\\
$w_{n,i} = w_{n,1}^{[i]}$ for $3 \leq i \leq n+1$,\\
\\
where, if $w_{n,j}$ is a word in $w_{n+2-i,1}$ and $w_{n+2-i,2}$,
 $w_{n,j}^{[i]}$ denotes
 the word obtained by switching $w_{n+2-i,1}$ and $w_{n+2-i,2}$.\\
\\
 This gives a well-defined sequence of words $w_{n,i}$
 for any positive integers $n, i \leq n+1$.\\
\\
The following formulas, which are just formal consequences of our
 notation, will be useful:
\begin{eqnarray*}
 w_{n,3} &=& w_{n,1}^{[3]}\\
    &=& W_n(w_{n-1,2}, w_{n-1,1})\\
    &=& (w_{n-1,2}^{p_n-1+q_n}w_{n-1,1}^{-q_n})^{k_n}
    w_{n-1,2}(w_{n-1,2}^{p_n-1+q_n}w_{n-1,1}^{-q_n})w_{n-1,2}^{-1} \,\,\,
   \textnormal{for } n \geq 2.
\end{eqnarray*}

\begin{eqnarray*}
 w_{n,i} &=& w_{n,1}^{[i]}\\
 &=& W_n(w_{n-1,1}^{[i-1]}, w_{n-1,2}^{[i-1]})\\
 &=& [(w_{n-1,1}^{[i-1]})^{p_n-1+q_n}(w_{n-1,2}^{[i-1]})^{-q_n}]^{k_n}\\
& & \hspace{1in} w_{n-1,1}^{[i-1]}[(w_{n-1,1}^{[i-1]})^{p_n-1+q_n}(w_{n-1,2}^{[i-1]})^{-q_n}]
 (w_{n-1,1}^{[i-1]})^{-1},\\
&& \textnormal{ for } n \geq 3 \textnormal{ and } 4 \leq i \leq n+1.
\end{eqnarray*}

Regardless of the choices of $a,b,p_i$ and $q_i$, we have:

\begin{pro} \label{sametrace}
Let $K$ be a field, and let $a,b \in SL_2(K)$.
Then $tr(w_{i,j}) =  tr(w_{i,k})$ for all $j,k \leq i+1$.
\end{pro}

We shall require the following lemma:

\begin{lem} \label{W}
$tr(W_n(x,y)) = tr(\bar{W}_n(x,y))$ for any $n$ and any\\
 $x,y \in K$.  
\end{lem}

\begin{proof}
Letting $U_n(x,y) = x^{k_n}y$, we have
\begin{eqnarray*}
tr(W_n(x,y)) &=& tr(U_n(x^{p_n-1+q_n}y^{-q_n},
                      x(x^{p_n-1+q_n}y^{-q_n})x^{-1}))\\
&=& tr(U_n(x(x^{p_n-1+q_n}y^{-q_n})x^{-1}, x^{p_n-1+q_n}y^{-q_n}))\\
& & \hspace{2in}   \textnormal{ (by Lemma \ref{Hor})}\\
&=& tr(\bar{W}_n(x,y)).
\end{eqnarray*}
\end{proof}

\begin{proof} (of Prop. \ref{sametrace})

We proceed by induction.
By Lemma \ref{W}, $tr(w_{1,1})= tr(w_{1,2})$.

Now, suppose that $tr(w_{n-1,i}) = tr(w_{n-1,j})$ for all $i,j \leq n$.

By Lemma \ref{W}, $tr(w_{n,1}) = tr(w_{n,2})$.
For $3 \leq i \leq n+1$, recall that
 $w_{n,1} = U(w_{n+2-i,1},w_{n+2-i,2})$ for some word $U$;  therefore
 we have:
\begin{eqnarray*}
tr(w_{n,1}) &=& tr(U(w_{n+2-i,1},w_{n+2-i,2})) \textnormal{ for some word }U\\
tr(w_{n,i}) &=& tr(w_{n,1}^{[i]})\\
&=& tr(U(w_{n+2-i,2},w_{n+2-i,1})).
\end{eqnarray*}
By the inductive hypothesis, $tr(w_{n+2-i,1})= tr(w_{n+2-i,2})$, so\\
 $tr(w_{n,1}) = tr(w_{n,i})$, by Lemma \ref{Hor}.
\end{proof}

\section{Algebraic representations}

The existence of maps from $\Gamma$ onto the groups $PSL_2(\mathbb{F}_{p_i})$
 depends on the existence of an algebraic representation of $\Gamma$,
 defined as follows:

 $\rho :\Gamma \rightarrow SL_2(\mathbb{C})$
 is an \textit{algebraic representation} if it is irreducible
 and its image is an infinite subgroup of $SL_2(\bar{\mathbb{Q}})$,
 where $\bar{\mathbb{Q}}$ is the algebraic closure of $\mathbb{Q}$.
  Note that this differs slightly from the definition given
 in [LR].

The purpose of this section is to show that any group $\Gamma$
 satisfying the hypotheses of Theorem \ref{trace} admits
 an algebraic representation.

\begin{lem} \label{algrep}
Let $\Gamma$ be a finitely generated group which admits an infinite
 irreducible representation $\rho_0 :\Gamma \rightarrow SL_2(\mathbb{C})$.
  Then $\Gamma$ admits an algebraic representation.
\end{lem}

\begin{proof} (of Lemma \ref{algrep})

Let $V_0(\Gamma)$ be the irreducible component of $V(\Gamma)$
 containing the character $\chi_{\rho_0}$ of the
 representation $\rho_0$ (recall the definition of $V(\Gamma)$ in Section 2).

 If $dim(V_0(\Gamma)) = 0$,
 then it is a well known fact that
 the coordinates of $\chi_{\rho_0}$ must be algebraic; in other
 words the image of $\chi_{\rho_0}$ must lie in $\bar{\mathbb{Q}}$.
  Then it follows from [Ba] that $\rho_0$
 is conjugate to an algebraic representation.

 So suppose $dim(V_0(\Gamma)) > 0$.  For $\gamma \in \Gamma$,
 define the function $\tau_{\gamma} : V_0(\Gamma) \rightarrow \mathbb{C}$
 by $\tau_{\gamma}(\chi_{\rho}) = \chi_{\rho}(\gamma) = tr \rho(\gamma) $.
 Recall, by [H] or [CS], that a character is determined by the
 values it takes on a finite set of elements
 $\gamma_1, ..., \gamma_n \in \Gamma$.
 Then since $dim(V_0(\Gamma)) > 0$, there
 is some $\gamma_i$ for which $\tau_{\gamma_i}$ is non-constant;
 let us assume it is $\gamma_1$.
 Since $\tau_{\gamma_1}$ is a non-constant polynomial map, it is 
 surjective.
 Now, all the characters in some neighborhood $U$ of $\chi_{\rho_0}$
 will correspond to infinite irreducible representations.
 $\tau_{\gamma_1}(U)$ is an open set in $\mathbb{C}$,
 and therefore contains an algebraic number $\alpha$.
 Let $\rho_1 \in \tau_{\gamma_1}^{-1}(\alpha)$; $\rho_1$ is infinite
 and irreducible, because it is in $U$.

 Suppose $\alpha$ has a minimal polynomial $f$ with coefficients
 in $\mathbb{Z}$.  Consider the algebraic subset
 $A_1(\Gamma) \subset V_0(\Gamma)$ obtained by adding the
 polynomial condition $f(\tau_{\gamma_1}(\chi_{\rho})) = 0$.
 $A_1(\Gamma)$ is non-empty, since it contains $\rho_1$.
  Let $V_1(\Gamma)$ denote the irreducible component of $A_1(\Gamma)$
 containing $\chi_{\rho_1}$.  Note that any character
 in $V_1(\Gamma)$ will map $\gamma_1$ to an algebraic number
 (in fact to a root of $f$).
 
  If $dim(V_1(\Gamma)) = 0$,
 then, as above, we have that $\rho_1$ is conjugate to an
 algebraic representation, and we are done.
  If $dim(V_1(\Gamma)) > 0$, then we can assume
 that $\tau_{\gamma_2}$ is non-constant on $V_1$,
 and then, as above, we can find an infinite irreducible
  representation $\rho_2 \in V_1$
 for which $\chi_{\rho_2}$ is algebraic.  Then we form the set $V_2(\Gamma)$,
 and so on.  Eventually, the process will terminate, when\\
\\
i. $dim(V_m(\Gamma)) = 0$ for some $m$,
 in which case we will get an algebraic representation, or\\
\\
ii. we have found an infinite irreducible representation $\rho$ such that
 $\chi_{\rho}(\gamma_i)$ is algebraic for all $i \leq n$.  Since $\chi_{\rho}$
 is a polynomial in the $\chi_{\rho}(\gamma_i)$'s, then
 $\chi_{\rho}(\gamma)$ is algebraic for all $\gamma \in \Gamma$,
 and therefore by [Ba], $\rho$ is conjugate to an algebraic representation.
\end{proof}

\section{Lemmas from group theory and number theory}

In this section we prove some lemmas which will be useful later,
 and we review the construction of the homomorphisms
 of $\Gamma$ into $PSL_2(\mathbb{C})$ alluded to in Section 3.
  For a more complete treatment of this construction, see [LR].
 For background on algebraic number theory, see [N].

Let $\rho$ be an algebraic representation of $\Gamma$ (see Section 5), and
 let $\mathbb{Q}(tr \rho(\Gamma))
       = \mathbb{Q}(\{tr\, \rho(\gamma) | \gamma \in \Gamma \})$.
  Since $\rho$ is algebraic,
 this is a finite extension of $\mathbb{Q}$.   It will be more
 convenient to work with the Galois closure, denoted
 $\overline{\mathbb{Q}(tr \rho \Gamma)}$.
 This is also a finite extension, of degree $N$, say.
 Let $\cal{O}$ denote the ring of integers of
 $\overline{\mathbb{Q}(tr \rho(\Gamma))}$.
  It follows from 
 the general theory of linear groups that for all but finitely
 many primes $p \in \mathbb{Z}$, there is a homomorphism
 $\phi_p : \Gamma \rightarrow PSL_2(\mathbb{F}_{p^i})$,
 where $\mathbb{F}_{p^i}$ is the residue field of a prime ideal
 $P \subset \cal{O}$ lying over $p$.  In particular,
 if $p$ \textit{splits completely} in
 $\overline{\mathbb{Q}(tr \rho(\Gamma))}$--
 i.e. factors into $N$ distinct prime ideals in $\cal{O}$--
 then $\phi_p$ maps into $PSL_2(\mathbb{F}_p)$.

    Since we shall require maps into groups
 of the form $PSL_2(\mathbb{F}_p)$, it will be useful to know
 how many primes in $\mathbb{Z}$ split completely in
 $\overline{\mathbb{Q}(tr \rho(\Gamma))}$.
  To give a precise answer requires the notion of natural density.

Let $A$ be a set of primes in $\mathbb{Z}$.  $A$ is said to have
 \textit{natural density} $\delta$ if
\begin{equation*}
\lim_{t \rightarrow \infty} 
\left( \frac{\# \textnormal{ of primes in } A < t}{\#
 \textnormal{ of primes in } \mathbb{Z} < t } \right) = \delta.
\end{equation*}

 We may now state a simple version of the Tchebotarev Density Theorem
 (see [L] p. 128).

\begin{thm} \label{Tcheb} (Tchebotarev density)
Let $K$ be a Galois extension of $\mathbb{Q}$ of degree $N$, and let
 $P = \{$ primes in $\mathbb{Z}$ which split completely in $K \}$.  Then
 the natural density of $P$ is $1/N$.  In particular, $P$ is infinite.
\end{thm}

In fact, if $p$ splits completely, then
 $\phi_p$ will almost always \textit{surject} onto
 $PSL_2(\mathbb{F}_p)$.  The following theorem follows directly
 from the proof of Theorem 1.2 in [LR]:

\begin{thm} \label{surj} \textnormal{(Long-Reid)}
 Let $\Gamma$ be a finitely generated group
 which admits an algebraic representation $\rho$.
 Let $\cal{P} = \{$ primes in $\mathbb{Z}$ which split completely
 in $\overline{\mathbb{Q}(tr \rho(\Gamma))}\}$, so that, for all but
 finitely many\\
 $p \in P$, the map
 $\phi_p: \Gamma \rightarrow PSL_2(\mathbb{F}_p)$ exists (see above).
  Then for all but finitely many $p \in \cal{P}$, $\phi_p $ is a surjection.
\end{thm}

We will in fact require $\Gamma$ to surject onto a product 
 of finite linear groups, prompting the following group theoretic
 lemma:

\begin{lem} \label{prod}
 Let $\Gamma$ be a group, and supppose $\Gamma$
 surjects onto a sequence $G_1, ..., G_n$ of distinct, finite simple
 groups.  Then $\Gamma$ surjects onto the direct
 product $\Pi_{i=1}^n G_i$.
\end{lem}

\begin{proof}
We are given surjections $\phi_i: \Gamma \rightarrow G_i$.\\
 Let $\phi_1 \times ... \times \phi_n : \Gamma \rightarrow \Pi_{i=1}^n G_i$
 denote the natural map induced by the $\phi_i$'s.
 Let $H$ denote the image of $\Gamma$ under
 $\phi_1 \times ... \times \phi_n$.

Let
\begin{equation*}
(*) \,\,\,\,\,
 1 = N_k \triangleleft N_{k-1} \triangleleft ...
 \triangleleft N_1 \triangleleft N_0 = H
\end{equation*}
be a chief series for $H$-- i.e. each $N_i$ is normal
 in $N_{i-1}$, and each quotient $N_{i-1}/N_i$ is simple.
The Jordan-Holder Theorem (see [I], p. 132) guarantees that such
 a series exists, and that the quotient groups
 are unique up to re-ordering.

Let $\pi_i:H \rightarrow G_i$ be the natural projection map.
 Since $\phi_i$ surjects $\Gamma$ onto $G_i$,
  $\pi_i$ surjects $H$ onto $G_i$.
Since $H/ ker \pi_i \cong G_i$ is simple,
it follows that $H$ has another chief series in which $ker \pi_i$
 is the first term (again see [I]).
Therefore for each $i$,
$G_i$ must appear as one of the quotients in the series (*).
Therefore
\begin{eqnarray*}
|H| &\geq& \Pi_{i=0}^{k-1} |N_i/N_{i+1}| \leq \Pi_{i=1}^n |G_i|\\
 &=& |\Pi_{i=1}^n G_i| 
\end{eqnarray*}
Therefore $H = \Pi_{i=1}^n G_i$, and so $\phi_1 \times ... \times \phi_n$
 is onto, proving the lemma.
\end{proof}

\begin{lem} \label{density}
Let $K$ be a Galois extension of degree $N$ over $\mathbb{Q}$,
 and let $\cal{P} = \{$ \textnormal{primes in } 
 $\mathbb{Z}$ \textnormal{ which split completely in} K $\}$. Let\\
 $p_1, ..., p_{n-1} \in \cal{P}$ with the property that
 $p_1 > 8N$ and $p_i > 2p_{i-1}$ for all $1 < i < n$.
Then, given any two integers $r$ and $s$, there exist infinitely many
 primes in $\cal{P}$ which
  are not congruent to $r$ or $s$ mod $p_{i}$ for any $i < n$.
\end{lem}

\begin{proof}
By Theorem \ref{Tcheb}, $\cal{P}$
 has natural density $1/N$ in $\mathbb{Z}$.

It follows from a more general version of Theorem \ref{Tcheb} (see [L], p.128)
 that if $p_i$ does not divide $r$ or $s$, then
the set of primes in $\mathbb{Z}$ which are congruent to $r$ mod
 $p_i$ has natural density $1/(p_i - 1)$, as does the set congruent to
 $s$ mod $p_i$.

Therefore, for large integers $m$ we have:

\begin{eqnarray*}
& &\frac{|\textnormal{primes }p < m : p \in P \textnormal{ and }
 p \not\equiv r
 \textnormal{ or } s \,\, (\textnormal{mod }p_1, ...\textnormal{ or }p_{n-1})|}
 {|\textnormal{primes} < m|}
\\
\\
&=& \frac{ |\textnormal{primes }p < m: p \in P|}{|\textnormal{ primes }p < m|}\\
& &  \hspace{.5in}- \frac{|\textnormal{primes }p < m : p \equiv r\textnormal{ or }s \,\,
 (\textnormal{mod } p_1, ...,\textnormal{ or } p_n)|}
{|\textnormal{ primes }p < m|}\\
& \geq& \frac{1}{|\textnormal{primes }p < m|}
[ |\textnormal{primes }p < m: p \in P|\\
& & - |\textnormal{primes }p < m: p \equiv r (\textnormal{mod }p_1)| - ...\\
& & - |\textnormal{primes }p < m: p \equiv r (\textnormal{mod }p_{n-1})|\\
& & - |\textnormal{primes }p < m: p \equiv s (\textnormal{mod }p_1)|\\
& & - ... - |\textnormal{primes }p < m: p \equiv s (\textnormal{mod }p_{n-1})|] \\
&=& \frac{1}{N} - \frac{1}{p_1-1}  - ... - \frac{1}{p_{n-1}-1}
 - \frac{1}{p_1-1} - ... - \frac{1}{p_{n-1}-1} \pm \epsilon, \\
& &\hspace{1cm}
\textnormal{ for some small $\epsilon$,  by the above density statements.}\\
&=& \frac{1}{N} - \frac{2}{p_1-1} - ... - \frac{2}{p_{n-1}-1} \pm \epsilon \\
&\geq& \frac{1}{N} - \frac{2}{8N} - ... - \frac{2}{2^{n+3}N}
     \pm \epsilon, \,\,\,\,\,\,\,\textnormal{ because }p_1 > 8N, \, 
   p_i > 2p_{i-1} \\
&>& \frac{1}{N} - \frac{1}{2N} \pm \epsilon \\
&=& \frac{1}{2N} \pm \epsilon. 
\end{eqnarray*}

So the ratio is bounded away from zero, and
therefore there must be infinitely many primes in $\cal{P}$ which
 are not congruent to either $r$ or $s$ mod $p_i$ for any $i < n$.
\end{proof} 

\section{Proving nonconjugacy}

In Section 4 we showed how to construct elements
 $w_{n,i} \in \Gamma$ in the same trace class.  Now we shall prove
 that these elements are pairwise non-conjugate.

It will be convenient to have the words written out in explicit form
here:
\\
\\
$w_{1,1} = (a^{p_1-1+q_1}b^{-q_1})^{k_1}a(a^{p_1-1+q_1}b^{-q_1})a^{-1}$\\
\\
$w_{1,2} = a(a^{p_1-1+q_1}b^{-q_1})^{k_1}a^{-1}(a^{p_1-1+q_1}b^{-q_1})$\\
\\
\\
$w_{2,1} = (w_{1,1}^{p_2 -1+ q_2}w_{1,2}^{-q_2})^{k_2}$
          $w_{1,1}(w_{1,1}^{p_2-1+q_2}w_{1,2}^{-q_2})w_{1,1}^{-1}$\\
\\
$w_{2,2} =w_{1,1}(w_{1,1}^{p_2 -1+ q_2}w_{1,2}^{-q_2})^{k_2}w_{1,1}^{-1}$
          $(w_{1,1}^{p_2-1+q_2}w_{1,2}^{-q_2})$\\
\\
$w_{2,3} = w_{2,1}^{[3]}  =(w_{1,2}^{p_2-1+q_2}w_{1,1}^{-q_2})^{k_2}$
          $w_{1,2}(w_{1,2}^{p_2-1+q_2}w_{1,1}^{-q_2})w_{1,2}^{-1}$\\        
\\
\\
$w_{n,1} = (w_{n-1,1}^{p_n-1 + q_n}w_{n-1,2}^{-q_n})^{k_n}$
     $w_{n-1,1}(w_{n-1,1}^{p_n-1 + q_n}w_{n-1,2}^{-q_n})w_{n-1,1}^{-1}$\\
\\
$w_{n,2}= w_{n-1,1}(w_{n-1,1}^{p_n-1+q_n}w_{n-1,2}^{-q_n})^{k_n}$
     $ w_{n-1,1}^{-1}(w_{n-1,1}^{p_n-1+q_n}w_{n-1,2}^{-q_n})$\\
\\
$w_{n,3} = w_{n,1}^{[3]} = W_n(w_{n-2,2}, w_{n-1,1})$\\
\indent $= (w_{n-1,2}^{p_n-1+q_n}w_{n-1,1}^{-q_n})^{k_n}
    w_{n-1,2}(w_{n-1,2}^{p_n-1+q_n}w_{n-1,1}^{-q_n})w_{n-1,2}^{-1}$,
   for $n \geq 2$\\
\\
$w_{n,i} = w_{n,1}^{[i]}$\\
\\       
$\indent = [(w_{n-1,1}^{[i-1]})^{p_n-1+q_n}(w_{n-1,2}^{[i-1]})^{-q_n}]^{k_n}$\\
$\indent\indent\indent\indent w_{n-1,1}^{[i-1]}[(w_{n-1,1}^{[i-1]})^{p_n-1+q_n}((w_{n-1,2}^{[i-1]})^{-q_n})]
        (w_{n-1,1}^{[i-1]})^{-1},$\\
\indent\indent\indent for $n \geq 3 \textnormal{ and } 4 \leq i \leq n+1$.\\
\\
 Let $\Gamma$ be as in the statement of Theorem \ref{trace},
 and let $\rho: \Gamma \rightarrow SL_2(\mathbb{C})$ be the algebraic
 representation guaranteed by Lemma \ref{algrep}.
 Let $\mathbb{Q}(tr \rho(\Gamma))$ be as in Section 6,
 let $\overline{\mathbb{Q}(tr \rho(\Gamma))}$ denote the Galois closure
 of $\mathbb{Q}(tr \rho(\Gamma))$, and let
 $N = |\overline{\mathbb{Q}(tr \rho(\Gamma))}:\mathbb{Q}|$.
 For primes $p$ which split completely in
 $\overline{\mathbb{Q}(tr \rho(\Gamma))}$, let the map 
 $\phi_p : \Gamma \rightarrow PSL_2(\mathbb{F}_p)$ be as discussed
 in Section 6.

 Let $\cal{Q} = \{ $ primes $p \in \mathbb{Z}$ for which
 $\phi_p$ is a surjection \}.
 Note that $\cal{Q}$ is infinite by Theorem \ref{surj}.
 The proof of Theorem \ref{trace} follows from:

\begin{pro} \label{nonconj}
The words $w_{n,i}$ can be chosen such that for every $n$
 there exist $a_n,b_n \in \Gamma$
 for which $w_{n,i}(a_n,b_n)$ is not
 conjugate to $w_{n,j}(a_n,b_n)$ or $w_{n,j}(a_n,b_n)^{-1}$
 whenever $i \neq j$.
 Indeed, the words $w_{n,i}$ and elements $a_n,b_n$
 can be chosen such that the following properties hold, for all $n$:\\
\\
P1(n): The primes $\{ p_1, p_2, ... p_n\}$ are in $\cal{Q}$,
 $p_1 > 8N$, and $p_i > 2p_{i-1}$ for all
 $i \leq n$.\\
\\
P2(n): For all $i \leq n$,
 $\phi_{p_i}(a_n) = [{2 \atop 0}{1 \atop 1/2}]$
 and $\phi_{p_i}(b_n) = [{2 \atop 0}{0 \atop 1/2}]$.\\
\\
P3(n): For all $i \leq n$ and $j \leq n+1-i$,
 $\phi_{p_i}(w_{n,j}(a_n,b_n)) = id$ and
 $\phi_{p_i}(w_{n,n+2-i}(a_n,b_n)) = [{1 \atop 0} {x \atop 1}] \neq id$.
\end{pro}

 We say that $p_i$ ``distinguishes'' $w_{n,j}(a_n,b_n)$
 from $w_{n,n+2-i}(a_n,b_n)$.

For example, taking $n = 3$, we have that $p_1$ distinguishes
 $w_{3,4}$ from $w_{3,3}, w_{3,2}$ and $w_{3,1}$;
 $p_2$ distinguishes $w_{3,3}$ from $w_{3,2}$ and $w_{3,1}$;
 and $p_3$ distinguishes $w_{3,2}$ from $w_{3,1}$.

\begin{proof} (of Prop. \ref{nonconj})
The proof is by induction. We begin by distinguishing $w_{1,1}$ from $w_{1,2}$.
We must first pick the prime $p_1$ and the integers $q_1$ and $k_1$
 which define $w_{1,1}$ and $w_{1,2}$.

 Let $q_1 = 1$,
 let $p_1 > max\{ 8N, 30 \}$ be a prime in $\cal{Q}$,
 and let $k_1 = p_1 - 4$.  $p_1$, $q_1$ and $k_1$ are now fixed, and will not
 change for the remainder of the proof.  Note that P1(1) is immediately
 satisfied.

 Let $a_1 \in \phi_{p_1}^{-1}([{2 \atop 0}{1 \atop 1/2}])$
 and $b_1 \in \phi_{p_1}^{-1}([{2 \atop 0} {0 \atop 1/2}])$,
 so P2(1) is satisfied.
  
Note that the images of $a_1$ and $b_1$ under $\phi_{p_1}$
 generate a semi-direct
 product of cyclic groups, that words of the form
 $[{1 \atop 0}{x \atop 1}]$ ($x \neq 0$) have order $p_1$,
 and that words of the form $[{z \atop 0} {x \atop z^{-1}}]$
 ($x\neq 0$, $z \neq 0,1$) have order dividing $p_1-1$.  This all follows
 from the structure theory of the groups $PSL_2(\mathbb{F}_p)$,
 and can be checked by explicit computation.

 We have:
\begin{eqnarray*}
\phi_{p_1}(w_{1,1}(a,b)) &=&
 \left[ \left( {2 \atop 0}\,\,{1 \atop 1/2} \right)^{p_1}
 \left( {1/2 \atop 0}\,\,{0 \atop 2} \right) \right]^{p_1-4}\\
& &\hspace{.25in}  \left( {2 \atop 0}\,\, {1 \atop 1/2} \right)
\left[ \left( {2 \atop 0}\,\,{1 \atop 1/2} \right)^{p_1}
 \left( {1/2 \atop 0}\,\,{0 \atop 2} \right) \right]
 \left( {1/2 \atop 0}\,\, {-1 \atop 2} \right)\\
 &=& \left[ \left( {2 \atop 0}\,\,{1 \atop 1/2} \right)
 \left( {1/2 \atop 0}\,\,{0 \atop 2} \right) \right]^{p_1-4}\\
& &\hspace{.25in}  \left( {2 \atop 0}\,\, {1 \atop 1/2} \right)
\left[ \left( {2 \atop 0}\,\,{1 \atop 1/2} \right)
 \left( {1/2 \atop 0}\,\,{0 \atop 2} \right) \right]
 \left( {1/2 \atop 0}\,\, {-1 \atop 2} \right),\\
&&\textnormal{since } \left( {2 \atop 0}{1 \atop 1/2} \right)
\textnormal{ has order dividing }p_1 -1.
\end{eqnarray*}
\begin{eqnarray*}
\indent \indent \indent &=& \left( {1 \atop 0}\,\,{2 \atop 1} \right)^{p_1-4}
    \left( {2 \atop 0}\,\,{1 \atop 1/2} \right)
    \left( {1 \atop 0}\,\,{2 \atop 1} \right)
    \left( {1/2 \atop 0}\,\,{-1 \atop 2} \right)\\
&=& \left( {1 \atop 0}\,\,{2p_1 - 8 \atop 1} \right)
    \left( {1 \atop 0}{8 \atop 1} \right)\\
&=& \left( {1 \atop 0}\,\,{2p_1 \atop 1} \right)\\
&=& \left({1 \atop 0}{0 \atop 1} \right).
\end{eqnarray*}

However,
\begin{eqnarray*}
\phi_{p_1}(w_{1,2}(a,b)) &=&
 \left( {2 \atop 0}\,\, {1 \atop 1/2} \right) 
 \left[ \left( {2 \atop 0}\,\,{1 \atop 1/2} \right)^{p_1}
 \left( {1/2 \atop 0}\,\,{0 \atop 2} \right) \right]^{p_1-4}\left( {1/2 \atop 0}\,\, {-1 \atop 2} \right)\\  
& & \hspace{.5in} \left[ \left( {2 \atop 0}\,\,{1 \atop 1/2} \right)^{p_1}
 \left( {1/2 \atop 0}\,\,{0 \atop 2} \right) \right]\\
&=&
\left( {2 \atop 0} \,\,{1 \atop 1/2} \right) 
 \left[ \left( {2 \atop 0}\,\,{1 \atop 1/2} \right)
 \left( {1/2 \atop 0}\,\,{0 \atop 2} \right) \right]^{p_1-4} \left( {1/2 \atop 0}\,\, {-1 \atop 2} \right)\\  
& &\hspace{.5in} \left[ \left( {2 \atop 0}\,\,{1 \atop 1/2} \right)
 \left( {1/2 \atop 0}\,\,{0 \atop 2} \right) \right]\\
&=&
\left( {2 \atop 0}\,\,{1 \atop 1/2} \right)
\left( {1 \atop 0}\,\,{2(p_1 - 4) \atop 1} \right)
\left( {1/2 \atop 0} \,\,{-1 \atop 2} \right)
\left( {1 \atop 0}{2 \atop 1} \right)\\
&=&
\left( {1 \atop 0}\,\, {8(p_1 - 4) \atop 1} \right)
\left( {1 \atop 0} {2 \atop 1} \right)\\
&=&
\left( {1 \atop 0} \,\,{-30 \atop 1} \right)\\
&\neq& \left( {1 \atop 0} {0 \atop 1} \right)
 \textnormal{since $p_1 > 30$.}
\end{eqnarray*}
Therefore $p_1$ distinguishes $w_{1,1}(a_1,b_1)$ from $w_{1,2}(a_1,b_1)$,
 and P3(1) is satisfied.

Now, suppose that we have
 picked $a_{n-1},b_{n-1} \in \Gamma$ and, for $i \leq n-1$,
 integers $k_i$, $p_i$, $q_i$ with corresponding words $w_{i,j}$,
 so that the following is true:\\
\\
P1(n-1): The primes $\{p_1, ..., p_{n-1} \}$ are in $\cal{Q}$,
 and $p_i > 2p_{i-1}$ for all $i \leq n-1$.\\
\\ 
P2(n-1): For all 
 $i \leq n-1$, $\phi_{p_i}(a_{n-1})=[{2 \atop 0} {1 \atop 1/2}]$
 and $\phi_{p_i}(b_{n-1}) = [{2 \atop 0} {0 \atop 1/2}]$.\\
\\ 
P3(n-1): For all $i \leq n-1$ and $j \leq n-i$,
 $\phi_{p_i}(w_{n-1,j}(a_{n-1},b_{n-1})) = id$ and
 $\phi_{p_i}(w_{n-1,n+1-i}(a_{n-1},b_{n-1}))=
 [{1 \atop 0} {x \atop 1}] \neq id$.\\
\\

 We shall show how to pick $k_n$, $q_n$, $p_n$, $a_n$ and $b_n$
 so that Properties P1(n), P2(n) and P3(n) are satisfied.\\
\\
Since $\cal{Q}$ is infinite by Theorem \ref{surj},
Property P1(n) may be satisfed simply by taking $p_n$ to be large enough.\\
\\
 Observe that, by
 Lemma \ref{prod},
 $\Gamma$ surjects onto $\Pi_{i=1}^nPSL_2(\mathbb{F}_{p_i})$.  Therefore we can
 satisfy Property P2(n) by picking
 $a_n \in \bigcap_{i=1}^n \phi_{p_i}^{-1}([{2 \atop 0}{1 \atop 1/2}])$
 and $b_n \in \bigcap_{i=1}^n \phi_{p_i}^{-1}([{2 \atop 0} {0 \atop 1/2}])$.\\
\\
To satisfy Property P3(n),
 we must show that, for all $i \leq n$ and $j \leq n+1-i$,
 $\phi_{p_i}(w_{n,j}) = id$ and
 $\phi_{p_i}(w_{n,n+2-i}) = [{1 \atop 0} {x \atop 1}] \neq id$.

Let $q_n = \Pi_{j=1}^{n-1}p_j(p_j-1)$. 

We will break the proof up into two cases; first we will see what happens
 when we reduce by the primes $p_i$, $i < n$, which we have already picked,
 and then we show how to pick $p_n$.\\
\\
\textit{Case 1}: $i < n$.\\
\\
$\phi_{p_i}(a_n)=[{2 \atop 0}{1 \atop 1/2}]$
 and $\phi_{p_i}(b_n)=[{2 \atop 0}{0 \atop 1/2}]$
 are contained in the subgroup  $B \subset PSL_2(\mathbb{F}_{p_i})$
 consisting of matrices whose lower left entry is 0.
 The order of $B$ is $p_i(p_i-1)/2$, which divides $q_n$.
  Therefore, $\phi_{p_i}(U(a_n,b_n)^{q_n}) = id$, where $U$ is any word
 on two letters.
 So we have, for $j \leq n+1-i$:\\
\\
\textit{Case 1a:} $j = 1$
\begin{eqnarray*}
 \phi_{p_i}(w_{n,1})
&=& \phi_{p_i}([w_{n-1,1}^{p_n-1 + q_n}w_{n-1,2}^{-q_n}]^{k_n}
    w_{n-1,1}[w_{n-1,1}^{p_n-1 + q_n}w_{n-1,2}^{-q_n}]w_{n-1,1}^{-1}\\
&=& \phi_{p_i}(w_{n-1,1}^{(p_n-1)(k_n+1)})\\
&=& id, \,\,\,\,\,\, \textnormal{by Property P3(n-1)}.
\end{eqnarray*}
\textit{Case 1b:} $j = 2$
\begin{eqnarray*}
\phi_{p_i}(w_{n,2})
 &=&\phi_{p_i}( w_{n-1,1}[w_{n-1,1}^{p_n-1+q_n}w_{n-1,2}^{-q_n}]^{k_n}
        w_{n-1,1}^{-1}[w_{n-1,1}^{p_n -1+ q_n} w_{n-1,2}^{-q_n}])\\
&=& \phi_{p_i}(w_{n-1,1}^{(p_n-1)(k_n+1)})\\
&=& id, \,\,\,\,\,\,\,\,\, \textnormal{by Property P3(n-1).}
\end{eqnarray*}
\textit{Case 1c:} $j = 3$
\begin{eqnarray*}
\phi_{p_i}(w_{n,3}) &=& \phi_{p_i}(w_{n,1}^{[3]})\\
&& \textnormal{(recall the definition of $w_{n,i}^{[j]}$ from Section 4)}\\
 &=& \phi_{p_i}(W_n(w_{n-2,2}, w_{n-1,1}))\\
 &=& \phi_{p_i}([w_{n-1,2}^{p_n-1+q_n}w_{n-1,1}^{-q_n}]^{k_n}
    w_{n-1,2}[w_{n-1,2}^{p_n-1+q_n}w_{n-1,1}^{-q_n}]w_{n-1,2}^{-1})\\
 &=& \phi_{p_i}(w_{n-1,2}^{(p_n-1)(k_n+1)})\\
 &=& id, \,\,\,\, \textnormal{by Property P3(n-1), since }j-1 = 2 \leq n-i.
 \end{eqnarray*}
\textit{Case 1d:} $4 \leq j \leq n+1-i$
\begin{eqnarray*}
\phi_{p_i}(w_{n,j}) &=& \phi_{p_i}(w_{n,1}^{[j]})\\
&=&\phi_{p_i}([(w_{n-1,1}^{[j-1]})^{p_n-1+q_n}(w_{n-1,2}^{[j-1]})^{-q_n}]^{k_n}\\
&&\hspace{.25in}w_{n-1,1}^{[j-1]}[(w_{n-1,1}^{[j-1]})^{p_n -1+ q_n}
                  (w_{n-1,2}^{[j-1]})^{-q_n}](w_{n-1,1}^{[j-1]})^{-1})\\
&=& \phi_{p_i}((w_{n-1,1}^{[j-1]})^{(p_n-1)k_n+p_n-1})\\
&=& \phi_{p_i}(w_{n-1,j-1}^{(p_n-1)(k_n+1)})\\
&=& id, \,\,\,\,\,\,\, \textnormal{by Property P3(n-1), since }j-1 \leq n-i.
\end{eqnarray*}
However,
\begin{eqnarray*}
\phi_{p_i}(w_{n,n+2-i}) &=& \phi_{p_i}(w_{n,1}^{[n+2-i]}), \,\,\,
 \textnormal{since }i\leq n-1, \textnormal{so }n+2-i \geq 3. \\
(\textnormal{for }i < n-1) &=& \phi_{p_i}([(w_{n-1,1}^{[n+1-i]})^{p_n-1+q_n}(w_{n-1,2}^{[n+1-i]})^{-q_n}]^{k_n}w_{n-1,1}^{[n+1-i]}\\
&&  [(w_{n-1,1}^{[n+1-i]})^{p_n-1+q_n}(w_{n-1,2}^{[n+1-i]})^{-q_n}]
  (w_{n-1,1}^{[n+1-i]})^{-1})\\
&=& \phi_{p_i}([w_{n-1,1}^{[n+1-i]}]^{(p_n-1)(k_n+1)})\\
&=& \phi_{p_i}(w_{n-1,n+1-i}^{(p_n-1)(k_n+1)})\\
&=& \left( {1 \atop 0} \,\, {(p_n-1)(k_n + 1)x \atop 1} \right)\\
&&\hspace{.5in} \textnormal{where } x \neq 0
 \textnormal{ by Property P3(n-1).}
\end{eqnarray*}
\begin{eqnarray*}
(\textnormal{for } i = n-1) &=& \phi_{p_{n-1}}(w_{n,1}^{[3]})\\
&=& \phi_{p_{n-1}}([w_{n-1,2}^{p_n-1+q_n}w_{n-1,1}^{-q_n}]^{k_n}\\
&&\hspace{.5in} w_{n-1,2}[w_{n-1,2}^{p_n-1+q_n}w_{n-1,1}^{-q_n}]w_{n-1,2}^{-1})\\
&=& \phi_{p_{n-1}}(w_{n-1,2}^{(p_n-1)(k_n+1)})\\
&=&  \,\,\,\,\left( {1 \atop 0} \,\, {(p_n-1)(k_n + 1)x \atop 1} \right),\\
& & \textnormal{ where } x \neq 0 \textnormal{ by Property P3(n-1)}.
\end{eqnarray*}

So to satisfy Property P3(n), we must pick
 $p_n \in \cal{Q}$ and $k_n$ such that:\\ 
\\
($I$) $p_i$ does not divide $k_n + 1$ or $p_n-1$ for $i < n$.\\
\\
Let $m$ be the sum of the exponents on $a_n$ and $b_n$
 in $w_{n-1,1}(a_n,b_n)$.  We shall also require:\\
\\
($II$) $p_n > max \{ 2^{4m}, q_n \}$.\\
\\
We set $k_n = p^n - 2^{2m}$.\\
\\
We claim that we can pick $p_n$ to satisfy properties ($I$) and ($II$).
  Indeed, ($I$) is equivalent to the statement\\
\\
$(I^{\prime})\,\,\, p_n \not \equiv 1$ or $2^{2m} - 1$ (mod $p_i$)
   for any $i < n$.\\
\\
Property P1(n-1) guarantees that $p_i > 2p_{i-1}$ for all $i < n$ and
 $p_1 > 8N$, so Lemma \ref{density}
 implies that $(I^{\prime})$ can be satisfied by infinitely
 many primes $p \in \cal{Q}$.  Therefore we can pick an arbitrarily large
 prime $p_n$
 to satisfy ($I$) and ($II$), and still satisfy Property P1(n).

This concludes the proof in Case 1.\\
\\
\textit{Case 2}: $i = n$\\
\\
We now show how to distinguish $w_{n,1}$ from $w_{n, 2}$.

 Since
 $\phi_{p_n}(a_n)=[{2 \atop 0}{1 \atop 1/2}]$
 and $\phi_{p_n}(b_n) = [{2 \atop 0}{0 \atop 1/2}]$, it is easy to
 see that  $\phi_{p_n}(w_{n-1,1})$
 will have the form $[{2^m \atop 0} {w \atop 2^{-m}}]$.  Note
 $2^m \neq 1$ (mod $p_n$), since $p_n > 2^{4m}$.
  Therefore, $\phi_{p_n}(w_{n-1,1}^{p_n-1}) = id$,
 and we have:
\begin{equation}
 \phi_{p_n}(w_{n,1}) = \phi_{p_n}((w_{n-1,1}^{p_n-1+q_n}
   w_{n-1,2}^{-q_n})^{k_n}
\end{equation}
\begin{equation*}
 \hspace{2in}w_{n-1,1}(w_{n-1,1}^{p_n-1+q_n}w_{n-1,2}^{-q_n})
   w_{n-1,1}^{-1})
\end{equation*}
\[ \indent \indent = \phi_{p_n}((w_{n-1,1}^{q_n}w_{n-1,2}^{-q_n})^{k_n}
   w_{n-1,1}(w_{n-1,1}^{q_n}w_{n-1,2}^{-q_n})w_{n-1,1}^{-1}). \] 

\begin{lem} \label{p_n}
If $p_n$ is chosen to be large enough, then
\[ \phi_{p_n}(w_{n-1,1}^{q_n}w_{n-1,2}^{-q_n})
  = \left( {1 \atop 0} {z \atop 1} \right), \]
 where $z \neq 0$.
\end{lem}

\begin{proof}

Recall that $p_{n-1}$ distinguishes $w_{n-1,1}$
 from $w_{n-1,2}$.  Therefore,
 $w_{n-1,1}(\bar{a},\bar{b}) \neq w_{n-1,2}(\bar{a},\bar{b})$
 for $\bar{a} = [{2 \atop 0} {1 \atop 1/2}]$,
 $\bar{b} = [{1 \atop 0} {1 \atop 1}]$
 $ \in PSL_2(\mathbb{F}_{p_{n-1}})$.

 Now let $\mathbb{Z}(1/2)$ denote the ring obtained by
 adjoining 1/2 to $\mathbb{Z}$, and let
 $\tilde{a}= [{2 \atop 0}{1 \atop 1/2}]$,
 $\tilde{b}=[{1 \atop 0}{1 \atop 1}] \in PSL_2(\mathbb{Z}(1/2))$.
 There is a well-defined reduction map
 $\theta_{p_{n-1}} : PSL_2(\mathbb{Z}(1/2)) \rightarrow
           PSL_2( \mathbb{F}_{p_{n-1}})$,
 since $p_{n-1} \neq 2$.
 Then
\begin{eqnarray*} \theta_{p_{n-1}}(w_{n-1,1}(\tilde{a},\tilde{b}))
      &=&  w_{n-1,1}(\bar{a},\bar{b})\\
      &\neq& w_{n-1,2}(\bar{a},\bar{b})
      = \theta_{p_{n-1}}(w_{n-1,2}(\tilde{a},\tilde{b})),
\end{eqnarray*}
so $w_{n-1,1}(\tilde{a}, \tilde{b}) \neq w_{n-1,2}(\tilde{a},\tilde{b})$
 in $PSL_2(\mathbb{Z}(1/2))$. 

 Using the definition of $w_{n-1,i}$ one may easily verify that
 the words have the following forms:
\begin{eqnarray*}
w_{n-1,1}(\tilde{a},\tilde{b}) =
    \left( {2^m \atop 0} {x \atop 2^{-m}} \right)\\
w_{n-1,2}(\tilde{a}, \tilde{b})) =
     \left( {2^m \atop 0} {y \atop 2^{-m}} \right),
\textnormal{ where } x \neq y.
\end{eqnarray*}

Then we compute:
\[ (w_{n-1,1}(\tilde{a},\tilde{b}))^{q_n} = 
                 \left( {2^m \atop 0} {x \atop 2^{-m}} \right)^{q_n}\]
\begin{eqnarray*}&=& \left( {2^{mq_n} \atop 0}
  {x(2^{m(q_n-1)}+2^{m(q_n-3)}+...+2^{m(3-q_n)}+2^{m(1-q_n)}) \atop 2^{-mq_n}}
   \right)\\
&\neq& \left( {2^{mq_n} \atop 0}
  {y(2^{m(q_n-1)}+2^{m(q_n-3)}+...+2^{m(3-q_n)}+2^{m(1-q_n)}) \atop 2^{-mq_n}}
  \right)\\
 &=& \left( {2^m \atop 0} {y \atop 2^{-m}} \right)^{q_n}\\
 &=& (w_{n-1,2}(\tilde{a},\tilde{b}))^{q_n}.
\end{eqnarray*}

 Then for a large enough prime $p_n$, these words are different
 mod $p_n$:
\begin{eqnarray*}
\phi_{p_n}(w_{n,1}(a,b))^{q_n} &=&
        \theta_{p_n}(w_{n,1}(\tilde{a},\tilde{b}))^{q_n}\\
 &=& \left( {v \atop 0}{x^{\prime} \atop 1/v} \right),
 \textnormal{ for some }v\\
\phi_{p_n}(w_{n,2}(a,b))^{q_n} &=&
        \theta_{p_n}(w_{n,2}(\tilde{a},\tilde{b}))^{q_n}\\
 &=& \left( {v \atop 0}{y^{\prime} \atop 1/v} \right)
\textnormal{, where } x^{\prime} \neq y^{\prime} \in \mathbb{F}_{p_n}.
\end{eqnarray*}

Then
\begin{eqnarray*}
\phi_{p_n}(w_{n,1}^{q_n}w_{n,2}^{-q_n})
&=& \left( {1 \atop 0}{v(x^{\prime} - y^{\prime}) \atop 1} \right)\\
&\neq& id,
\end{eqnarray*}
completing the proof of Lemma \ref{p_n}.
\end{proof}

Then, returning to Equation 2, we have
\begin{eqnarray*}
\phi_{p_n}(w_{n-1,1}^{q_n}w_{n-1,2}^{-q_n})
  = \left( {1 \atop 0} {z \atop 1} \right) \neq id.
\end{eqnarray*}
So,
\begin{eqnarray*}
\phi_{p_n}(w_{n,1}) &=& \left( {1 \atop 0} {k_nz \atop 1} \right)
  \left( {2^m \atop 0} \,\, {w \atop 2^{-m} } \right)
  \left( {1 \atop 0} \,\, {z \atop 1} \right)
  \left( {2^{-m} \atop 0}\,\, {-w \atop 2^m } \right) \textnormal{ for some }
 w.\\
&=& \left( {1 \atop 0} \,\, {k_nz \atop 1} \right)
  \left( {1 \atop 0} {2^{2m}z \atop 1} \right)\\
&=& \left( {1 \atop 0} \,\, {(k_n + 2^{2m})z \atop 1} \right)\\
&=& \left( {1 \atop 0} \,\, {(p_n - 2^{2m} + 2^{2m})z \atop 1} \right),
  \textnormal{ by our choice of }k_n.\\  
&=& \left( {1 \atop 0} { 0 \atop 1} \right).
\end{eqnarray*}

However:
\begin{eqnarray*}
\phi_{p_n}(w_{n,2}) &=&
  \phi_{p_n}(w_{n-1,1}(w_{n-1,1}^{p_n-1+q_n}w_{n-1,2}^{-q_n})^{k_n}
      w_{n-1,1}^{-1}(w_{n-1,1}^{p_n-1+q_n}w_{n-1,2}^{-q_n}))\\
&=& \phi_{p_n}(w_{n-1,1}(w_{n-1,1}^{q_n}w_{n-1,2}^{-q_n})^{k_n}
      w_{n-1,1}^{-1}(w_{n-1,1}^{q_n}w_{n-1,2}^{-q_n}))\\
&=& \left( {2^m \atop 0}\,\, {w \atop 2^{-m}} \right)
  \left( {1 \atop 0}\,\, {k_nz \atop 1} \right)
  \left( {2^{-m} \atop 0} \,\,{-w \atop 2^m} \right)
  \left( {1 \atop 0}\,\, {z \atop 1} \right)\\
&=& \left( {1 \atop 0}\,\, { 2^{2m}k_nz + z \atop 1} \right)\\
&=& \left( {1 \atop 0}\,\, { (2^{2m}(p_n - 2^{2m}) +1)z \atop 1} \right)\\
&=& \left( {1 \atop 0}\,\, { (-2^{4m} + 1)z \atop 1} \right).
\end{eqnarray*}
Since $p_n > 2^{4m}$, $-2^{4m} + 1 \neq 0$ (mod $p_n$).
 Also, $z \neq 0$,
 so $\phi_{p_n}(w_{n,2}) \neq 0$.  Hence
 Property P3(n) is satisfied,
 and we are done.

\end{proof}

\section{references}

\vspace{1pc}

[A] J. Anderson, ``A brief survey of the deformation theory
 of Kleinian groups'', preprint.

\vspace{1pc}

[Ba] H. Bass, Groups of integral representation type, Pacific J. Math.
 \textbf{86} (1980), 15-51.

\vspace{1pc}

[Br] R. Brooks, ``Circle packings and co-compact extensions
 of Kleinian groups'', Invent. Math. \textbf{86} (1986), 461-469.

\vspace{1pc}

[CS] M. Culler and P. B. Shalen, ``Varieties of group representations
 and splittings of 3-manifolds'', Ann. Math. \textbf{117} (1983),
 109-146.

\vspace{1pc}

[GR] D. Ginzburg and Z. Rudnick, ``Stable multiplicities in the
 length spectrum of Riemann surfaces'', Isr. J. Math. \textbf{104} (1998),
 129-144.
 
\vspace{1pc}

[H] R. D. Horowitz, ``Characters of free groups represented in the
 two-dimensional special linear group'', Comm. Pure Appl. Math \textbf{25}
 (1972), 635-649.

\vspace{1pc}

[I] I. M. Isaacs, Algebra: a Graduate Course, Brooks/Cole, 1994.

\vspace{1pc}

[L] R. L. Long, Algebraic Number Theory, Marcel Dekker, 1977.

\vspace{1pc}

[LR] D. D. Long and A.W. Reid, ``Simple quotients of Hyperbolic
 3-manifold groups'', Proc. Amer. Math. Soc. \textbf{126},
 (1998), 877-880.

\vspace{1pc}

[M] J. Marklof, ``On multiplicities in length spectra of arithmetic
 hyperbolic three-orbifolds'', Nonlinearity \textbf{9} (1996), 517-536.

\vspace{1pc}

[N] W. Narkiewicz, Algebraic Numbers, Polish Scientific Publishers,
 Warsaw, 1974.

\vspace{1pc}

[R] B. Randol, ``The length spectrum of a Riemann surface is always
 of unbounded multiplicity'', Proc. Amer. Math. Soc. \textbf{78} (1980),
 455-456.

\vspace{1pc}

[Sar] P. Sarnak, ``Arithmetic quantum chaos'', The Schur Lectures (1992),
 Israel Mathematical Conference Proceedings  Vol. 8, 1995.

\vspace{1pc}

[Sc] P. Scott, ``Compact submanifolds of 3-manifolds'', J. London Math. Soc.
 \textbf{7} (1973), 246-250.

\vspace{1pc}

[Su] M. Suzuki, Group Theory I, Springer-Verlag, 1982.

\vspace{1pc}

[T] W. Thurston, Geometry and Topology of Hyperbolic\\
 3-manifolds, mimeographed lecture notes, 1978.

\newpage

\noindent{Department of Mathematics}\\
University of Texas at Austin\\
Austin TX 78712\\
\textit{masters@math.utexas.edu}\\
\\
Current Address:\\
Department of Mathematics\\
Rice University\\
Houston TX 77005-1892\\
\textit{masters@math.rice.edu}
\end{document}